\documentstyle[11pt]{article}
\def\b#1{{\bf #1}}
\def\i#1{{\it #1}}

\begin{document}
\title{Generalized Bivariate Fibonacci Polynomials}
\author{Mario Catalani\\
Department of Economics, University of Torino\\ Via Po 53, 10124 Torino, Italy\\
mario.catalani@unito.it}
\date{}
\maketitle
\begin{abstract}
\small{We define generalized bivariate polynomials, from which specifying
initial conditions the bivariate Fibonacci and
Lucas polynomials are obtained. Using essentially a matrix approach we derive
identities and inequalities that in most cases generalize known results.}
\end{abstract}

\section{Antefacts}
The generalized bivariate Fibonacci polynomial may be defined as
$$H_n(x,\,y)=xH_{n-1}(x,\,y)+yH_{n-2}(x,\,y),\qquad  H_0(x,\,y)=a_0,\,
H_1(x,\,y)=a_1.$$
We assume $y\not= 0$, as well as $x^2+4y\not= 0$.

\noindent
If we set $a_0=0,\,a_1=1$ we obtain the bivariate Fibonacci polynomials
$F_n(x,\,y)$;
with $a_0=2,\,a_1=x$ we get the bivariate Lucas polynomials.
$L_n(x,\,y)$. The univariate polynomials are obtained setting $y=1$. The
characteristic polynomial is
$$f(t)=t^2-tx-y,$$
and the roots are (see for example \cite{gauthier1})
$$\alpha\equiv \alpha(x,\,y)={x+\sqrt{x^2+4y}\over 2},\qquad
\beta\equiv \beta(x,\,y)={x-\sqrt{x^2+4y}\over 2}.$$
Note that $\alpha +\beta=x$, $\alpha\beta=-y$. Furthermore $\alpha>\beta$;
if $y>0$ then $\alpha>0$ and $\beta <0$ and if $x>0$ then $\alpha>
\vert\beta\vert$, while if
$x<0$ then $\alpha<
\vert\beta\vert$;
if $y<0$ then $\alpha$ and $\beta$ have the same sign: if $x>0$ they are
both positive, while if
$x<0$ they are both negative.

\noindent
The ordinary generating function $g(t)$ is
$$g(t)={a_0+(a_1-a_0x)t\over 1-xt-yt^2},$$
which can be obtained through the usual trick: write
$$g(t)=H_0+H_1t+H_2t^2+\cdots $$
then obtain $xtg(t)$ and $yt^2g(t)$, subtract them from $g(t)$, use the
recurrence to get rid of all summands except the first two and impose the
initial conditions.

\noindent
Using the Rational Expansion Theorem (see
\cite{knuth}) we get the Binet's form
\begin{equation}
H_n(x,\,y)=(a_1-\beta a_0){\alpha^n\over \alpha -\beta}-
(a_1-\alpha a_0){\beta^n\over \alpha -\beta}.
\end{equation}
Hence
$$F_n(x,\,y)={\alpha^n\over \alpha -\beta}-{\beta^n\over \alpha -\beta},$$
$$L_n(x,\,y)=\alpha^n+\beta^n.$$
In the sequel we will write simply $H_n,\, F_n$ and $L_n$ instead of
$H_n(x,\,y),\, F_n(x,\,y)$ and
$L_n(x,\,y)$.

\noindent
The sequence can be extended to negative subscripts by defining
$$H_{-n}=-{x\over y}H_{-(n-1)}+{1\over y}H_{-(n-2)}.$$
In this case the characteristic polynomial is the reflected polynomial
(see \cite{knuth}) of $f(t)$: hence the roots are ${1\over\alpha}$ and
${1\over\beta}$.
Define the matrix
\begin{equation}
\label{eq:unouno}
\b{A}=\left [\begin{array}{cc}x&1\\y&0\end{array}\right ].
\end{equation}
Note that $\vert\b{A}\vert = -y$.

\noindent
By easy induction
\begin{equation}
\label{eq:matriceuno}
\b{A}^n=\left [\begin{array}{cc}F_{n+1}&F_n\\
yF_n&yF_{n-1}\end{array}\right ].
\end{equation}
From the above equation taking determinants
\begin{equation}
\label{eq:fibonacciuno}
F_{n-1}F_{n+1}-F_n^2=(-1)^ny^{n-1}.
\end{equation}
Also from $\b{A}^{2n}=\b{A}^n\b{A}^n$ it follows
\begin{equation}
F_{n+1}^2+yF_n^2=F_{2n+1}.
\end{equation}
The eigenvalues of \b{A} are $\alpha$ and $\beta$: from this
\begin{equation}
\label{eq:fibonaccidue}
{\rm tr}(\b{A}^n)=F_{n+1}+yF_{n-1}=\alpha^n+\beta^n=L_n.
\end{equation}
Also it is immediate to get
$$F_nL_n=F_{2n}.$$
\noindent
A wealth of inequalities can be obtained from Equation~\ref{eq:matriceuno}
using the classical matrix inequality on the extremes of quadratic forms
\begin{equation}
\beta^n(z_1^2+z_2^2)\le z_1^2F_{n+1}+z_1z_2(1+y)F_n+z_2^2yF_{n-1}
\le \alpha^n(z_1^2+z_2^2),
\end{equation}
for all real $z_1,\,z_2$.

\noindent
Now define the matrix
\begin{eqnarray*}
\b{C}&=&\left [\begin{array}{cc}a_1x+a_0y&a_1\\a_1y&a_0y\end{array}\right ]\\
&=&a_0y\b{I}+a_1\b{A}.
\end{eqnarray*}
Again by easy induction we get
\begin{equation}
\label{eq:matricegenerale}
\b{CA}^n
=\left [\begin{array}{cc}H_{n+2}&H_{n+1}\\
yH_{n+1}&yH_n\end{array}\right ].
\end{equation}
Taking $a_0=2,\,a_1=x$ we get the matrix
\begin{eqnarray*}
\b{B}&=&\left [\begin{array}{cc}x^2+2y&x\\xy&2y\end{array}\right ]\\
&=&2y\b{I}+x\b{A},
\end{eqnarray*}
and
\begin{equation}
\b{BA}^n
=\left [\begin{array}{cc}L_{n+2}&L_{n+1}\\
yL_{n+1}&yL_n\end{array}\right ].
\end{equation}
Because
$$\b{CA}^n=(a_0y\b{I}+a_1\b{A})\b{A}^n=a_0y\b{A}^n+a_1\b{A}^{n+1},$$
we obtain the identity
$$H_{n+1}=a_0yF_n+a_1F_{n+1},$$
which specializes to
$$L_{n+1}=2yF_n+xF_{n+1}.$$
Analogous results we get for $H_{-n}$ using the inverse of \b{A}
$$\b{A}^{-1}=\left [\begin{array}{cc}0&{1\over y}
\\1&{-{x\over y}}\end{array}\right ].$$
In particular we get
$$H_{-n}={(-1)^n\over y}\left [a_0F_{n+1}-a_1F_n\right ].$$
Because $\vert \b{C}\vert = a_0a_1xy+a_0^2y^2-a_1^2y$ taking determinants in
Equation~\ref{eq:matricegenerale} we get
\begin{equation}
\label{eq:lucasuno}
H_nH_{n+2}-H_{n+1}^2=(-1)^ny^{n-1}(a_0a_1xy+a_0^2y^2-a_1^2y),
\end{equation}
which can be termed the \i{Simpson} formula.

\noindent
Using identity 1.60 in \cite{gould}, that is
$$\sum_{k=0}^{\left\lfloor{n\over 2}\right\rfloor}
(-1)^k{n-k\choose k}(xy)^k(x+y)^{n-2k}={x^{n+1}-y^{n+1}\over x-y},$$
writing $\alpha$ for $x$ and $\beta$ for $y$ we have
\begin{eqnarray*}
\sum_{k=0}^{\left\lfloor{n\over 2}\right\rfloor}
(-1)^k{n-k\choose k}(\alpha\beta)^k(\alpha+\beta)^{n-2k}&=&
{\alpha^{n+1}-\beta^{n+1}\over \alpha-\beta}\\
&=&F_{n+1},
\end{eqnarray*}
which means
$$\sum_{k=0}^{\left\lfloor{n\over 2}\right\rfloor}
{n-k\choose k}y^kx^{n-2k}=
F_{n+1}.$$

\section{Inversion formulas}
Now we will establish the following formula for $\gamma$ complex and
${\rm i}^2=-1$
\begin{equation}
F_n(x,\,-y)=({\rm i}\gamma)^{-(n-1)}F_n({\rm i}\gamma x,\,\gamma^2y).
\end{equation}
We have
$$F_n(x,\,-y)={\alpha_1^n-\beta_1^n\over\alpha_1-\beta_1},$$
where
$$\alpha_1={x+\sqrt{x^2-4y}\over 2}.$$
Analogously for $\beta_1$. Also
$$F_n({\rm i}\gamma x,\,\gamma^2y)=
{\alpha_2^n-\beta_2^n\over\alpha_2-\beta_2},$$
where
\begin{eqnarray*}
\alpha_2&=&{{\rm i}\gamma x+\sqrt{-\gamma^2x^2+4\gamma^2y}\over 2}\\
&=&{{\rm i}\gamma x+{\rm i}\gamma\sqrt{x^2-4y}\over 2}\\
&=&{\rm i}\gamma{x+\sqrt{x^2-4y}\over 2}\\
&=&{\rm i}\gamma\alpha_1.
\end{eqnarray*}
Analogously for $\beta_2$.
Then
\begin{eqnarray*}
F_n({\rm i}\gamma x,\,\gamma^2y)&=&
{({\rm i}\gamma)^n\alpha_1^n-({\rm i}\gamma)^n\beta_1^n
\over{\rm i}\gamma\alpha_1-{\rm i}\gamma\beta_1}\\
&=&({\rm i}\gamma)^{n-1}F_n(x,\,-y),
\end{eqnarray*}
whence the desired formula.
Some consequences. If we set $\gamma={\rm i}$ we get
\begin{equation}
F_n(x,\,-y)=(-1)^{n-1}F_n(-x,\,-y).
\end{equation}
Also from
$$F_n(x,\, y)=({\rm i}\gamma)^{-(n-1)}
F_n({\rm i}\gamma x,\,-\gamma^2 y)$$
we get with $\gamma={1\over \sqrt{- y}}={1\over {\rm i}
\sqrt{y}}$
\begin{equation}
\label{eq:duedue}
F_n(x,\,y)=\left (\sqrt{ y}\right )^{n-1}F_n\left ({x\over
\sqrt{ y}},\,1\right ).
\end{equation}
In a similar fashion we obtain
\begin{equation}
L_n(x,\,-y)=({\rm i}\gamma)^{-n}L_n({\rm i}\gamma x,\,\gamma^2y),
\end{equation}
$$L_n(x,\,y)=\left (\sqrt{ y}\right )^nL_n\left ({x\over
\sqrt{ y}},\,1\right ),$$
and
$$L_n(-x,\,y)=(-1)^nL_n(x,\,y).$$
It follows
$$F_{2n+1}(-x,\,y)=F_{2n+1}(x,\,y),\quad
F_{2n}(-x,\,y)=-F_{2n}(x,\,y),$$
$$L_{2n+1}(-x,\,y)=-L_{2n+1}(x,\,y),\quad
L_{2n}(-x,\,y)=F_{2n}(x,\,y).$$

From the Binet's forms we get
$$F_n(ax,\,by)=a^{n-1}F_n\left (x,\,{b\over a^2}y\right )$$
$$L_n(ax,\,by)=a^nL_n\left (x,\,{b\over a^2}y\right )$$
from which we get
$$F_n(x,\,y)=x^{n-1}F_n\left (1,\,{y\over x^2}\right ),$$
$$L_n(x,\,y)=x^nL_n\left (1,\,{y\over x^2}\right ).$$

\section{A complex representation}
Consider the matrix
$$\b{B}={1\over \sqrt{x+2{\rm i}}}
\left [\begin{array}{cc}x+{\rm i}&1\\1&{\rm i}\end{array}\right ],$$
where ${\rm i}=\sqrt{-1}$. Then
$$\b{B}^2=\left [\begin{array}{cc}x &1\\1&0\end{array}\right ]=\b{A}.$$
This matrix is obtained by Equation~\ref{eq:unouno} setting $y=1$ and so
is the matrix generating the univariate Fibonacci polynomials.
Then
$$\b{A}^n=\left [\begin{array}{cc} F_{n+1}(x,\,1)&F_n(x,\,1)\\
F_n(x,\,1)&F_{n-1}(x,\,1)\end{array}\right ].$$
Using the
method exposed in \cite{mclaughlin} to $\b{B}^{2n}=\b{A}^n$ (which
consists in expressing the elements of the $n$-th power of a $2\times 2$
matrix in terms of the elements of the matrix)
we get
\begin{equation}
F_n(x,\,1)=\sum_{r=0}^{n-1}{2n-1-r\choose r}(x+2{\rm i})^{n-r-1}
(-{\rm i})^r.
\end{equation}
Using Equation~\ref{eq:duedue},
after replacing in the previous equation $x$ by ${x\over\sqrt{y}}$, we get
\begin{equation}
F_n(x,\,y)=\sum_{r=0}^{n-1}
{2n-1-r\choose r}(x+2\sqrt{y}{\rm i})^{n-r-1}(-\sqrt{y}{\rm i})^r.
\end{equation}

\section{Summations}
Since
$$\sum_{k=0}^n\alpha^k={1-\alpha^{n+1}\over 1-\alpha},$$
(analogously for $\beta$)
and $(1-\alpha)(1-\beta)=1-x-y$, we have
\begin{eqnarray*}
\sum_{k=0}^nH_k&=&{a_1-a_0\beta\over\alpha-\beta}\sum_{k=0}^n\alpha^k
-{a_1-a_0\alpha\over\alpha-\beta}\sum_{k=0}^n\beta^k\\
&=&{a_0(1-x)+a_1-yH_n-H_{n-1}\over 1-x-y},
\end{eqnarray*}
which can be rewritten as
$$(1-x-y)\sum_{k=0}^nH_k=
a_0(1-x)+a_1-yH_n-H_{n-1}.$$
Specializing we get
$$(1-x-y)\sum_{k=0}^nF_k=1-yF_n-F_{n+1},$$
$$(1-x-y)\sum_{k=0}^nL_k=2-x-yL_n-L_{n+1}.$$
Since
$$\sum_{k=0}^n\alpha^{2k}={1-\alpha^{2n+2}\over 1-\alpha^2},$$
(analogously for $\beta$),
$(1-\alpha^2)(1-\beta^2)=1-x^2-2y+y^2$, and $\alpha^2-\beta^2=x(\alpha-
\beta)$ we have
$$(1-2y-x^2+y^2)\sum_{k=0}^nF_{2k}=x+y^2F_{2n}-F_{2n+2},$$
$$(1-2y-x^2+y^2)\sum_{k=0}^nL_{2k}=2-x^2-2y+y^2L_{2n}-L_{2n+2}.$$
Just a little bit more involved is the following inequality. Consider
\begin{eqnarray*}
L_n^m&=&(\alpha^n+\beta^n)^m\\
&=&\sum_{k=0}^m{m\choose k}\alpha^{nk}\beta^{n(m-k)}.
\end{eqnarray*}
Now upon expanding the summation and grouping the first and the last
element, the second and the next-to-last element and so on we get the
following formulas depending on the parity of $m$
$$L_n^m=\sum_{k=0}^{\left\lfloor{m\over 2}\right\rfloor}
{m\choose k}(-y)^{kn}L_{n(m-2k)}\qquad {\rm if}\,\, \mbox{$m$ is odd,}$$
$$L_n^m=\sum_{k=0}^{{m\over 2}-1}
{m\choose k}(-y)^{kn}L_{n(m-2k)}+ {m\choose {m\over 2}}(-y)^{mn\over 2}
\qquad \mbox{if $m$ is even.}$$
If we define
$$\varepsilon_1=\left\lfloor{m\over 2}\right\rfloor+{1\over 2}
\left [-1+(-1)^{m+1}\right ], \qquad \varepsilon_2 =
{1\over 2}
\left [(-1)^m+1\right ],$$
we can group the two identities into the following one
$$L_n^m=\sum_{k=0}^{\varepsilon_1}
{m\choose k}(-y)^{kn}L_{n(m-2k)}+ \varepsilon_2{m\choose
{m\over 2}}(-y)^{mn\over 2}.$$
Other consequences are
$$\sum_{k=0}^n{n\choose k}x^ky^{n-k}F_k=F_{2n},$$
$$\sum_{k=0}^n{n\choose k}x^ky^{n-k}L_k=L_{2n}.$$
Indeed, for instance the second one, we have
\begin{eqnarray*}
\sum_{k=0}^n{n\choose k}x^ky^{n-k}L_k&=&
\sum_{k=0}^n{n\choose k}x^ky^{n-k}\alpha^k +
\sum_{k=0}^n{n\choose k}x^ky^{n-k}\beta^k\\
&=&(y+x\alpha)^n+(y+x\beta)^n,
\end{eqnarray*}
and $y+x\alpha=\alpha^2$, $y+x\beta=\beta^2$.

\section{Eigenvectors and eigenvalues}
Now $\forall\, k,\,n$
\begin{eqnarray*}
\b{A}^k\b{BA}^n&=&
\b{A}^k(y\b{I}+\b{A}^2)\b{A}^n\\
&=&y\b{A}^{n+k}+\b{A}^{n+k+2}\\
&=&(y\b{I}+\b{A}^2)\b{A}^{n+k}\\
&=&\b{BA}^{n+k}.
\end{eqnarray*}
This means that $\forall\, k,\,n$ the matrices $\b{A}^k$ and $\b{BA}^n$
commute: then (see \cite{mario}, \cite{gantmacher},
\cite{rao}) the matrices
possess the same eigenvectors.

\noindent
The eigenvectors (not normalized) of \b{A} are
$$\b{z}_1= \left [\begin{array}{c}{\alpha\over y}\\ 1\end{array}\right ],
\quad {\rm and }\quad
\b{z}_2= \left [\begin{array}{c}{\beta\over y}\\ 1\end{array}\right ].$$
For the eigenvalues of $\b{BA}^n$ since
$$\b{BA}^n = y\b{A}^n+\b{A}^{n+2}$$
we have
$$\left \{\begin{array}{c}\lambda_1=y\alpha^n+\alpha^{n+2},\\
\lambda_2=y\beta^n+\beta^{n+2}.
\end{array}\right .$$
Consider $(\b{BA}^n)^2=\b{BA}^n\b{BA}^n$. Since
$\b{A}^k$ and $\b{BA}^n$
commute, this means that we can write $(\b{BA}^n)^2=\b{B}\b{BA}^{2n}$, that
is
$$
\left [\begin{array}{cc}2y+x^2&x\\xy&2y\end{array}\right ]
\left [\begin{array}{cc}L_{2n+2}&L_{2n+1}\\
yL_{2n+1}&yL_{2n}\end{array}\right ].$$
Using repeatedly the recurrence for $L_n$ we get
\begin{equation}
L_{n+2}^2+yL_{n+1}^2=L_{2n+4}+yL_{2n+2}.
\end{equation}

\noindent
Another immediate consequence is the following. Since
$$\b{A}^n\b{BA}^n=
\left [\begin{array}{cc}F_{n+1}L_{n+2}+yF_nL_{n+1}&F_{n+1}L_{n+1}+
yF_nL_n\\
yF_nL_{n+2}+y^2F_{n-1}L_{n+1}&yF_nL_{n+1}+y^2F_{n-1}L_n\end{array}\right ],$$
and $\b{A}^n\b{BA}^n=\b{BA}^{2n}$ it follows equating the elements in
position $(2,\,2)$
\begin{equation}
L_{2n}=F_nL_{n+1}+yF_{n-1}L_n.
\end{equation}
Because of the \i{Schur inequality} (see \cite{mario}, \cite{marcus}),
that is: the sum of the squared elements is greater or equal to the
sum of the squared eigenvalues, we have the following inequality
\begin{eqnarray*}
F_{n+1}^2+(1+y^2)F_n^2+y^2F_{n-1}^2&\ge &\alpha^{2n}+\beta^{2n}\\
&=&L_{2n}.
\end{eqnarray*}
Note that the equality sign holds if $y=1$.

\noindent
In the same way we get
\begin{eqnarray*}
y^2L_n^2+(1+y^2)L_{n+1}^2+L_{n+2}^2&\ge&(y\alpha^n+\alpha^{n+2})^2
+(y\beta^n+\beta^{n+2})^2\\
&=&y^2L_{2n}+2yL_{2n+2}+L_{2n+4},
\end{eqnarray*}
and the equality sign holds if $y=1$.

\noindent
Consider
\begin{eqnarray*}
\b{B}^2\b{A}^n&=&
y^2\b{A}^n+\b{A}^{n+4}+2y\b{A}^{n+2}\\
&=&y\b{BA}^n+\b{BA}^{n+2}\\
&=&
\left [\begin{array}{cc}yL_{n+2}+L_{n+4}&yL_{n+1}+L_{n+3}\\
y^2L_{n+1}+yL_{n+3}&y^2L_n+yL_{n+2}\end{array}\right ].
\end{eqnarray*}
Taking the determinants we get
\begin{eqnarray*}
(-1)^ny^n(x^2y+4y^2)^2&=&y^3L_nL_{n+2}+y^2L_{n+2}^2+y^2L_nL_{n+4}+
yL_{n+2}L_{n+4}\\
&&\qquad -y^3L_{n+1}^2-yL_{n+3}^2-2y^2L_{n+1}L_{n+3}.
\end{eqnarray*}
Using repeatedly identity~\ref{eq:lucasuno} we obtain
\begin{equation}
L_nL_{n+4}-L_{n+1}L_{n+3}=(-1)^n(x^4y^n+5x^2y^{n+1}+4y^{n+2}).
\end{equation}

\end{document}